\documentclass[a4paper, 12pt, twoside]{article}
\usepackage[left=2.9cm,right=2.9cm, top=3.4cm,bottom=3.6cm,bindingoffset=0cm]{geometry}

\newenvironment{resetCounters}
{
\setcounter{section}{0}
\setcounter{subsection}{0}
\setcounter{subsubsection}{0}
\setcounter{paragraph}{0}
\setcounter{subparagraph}{0}
\setcounter{equation}{0}
\setcounter{figure}{0}
\setcounter{table}{0}
\setcounter{footnote}{0}
\setcounter{theorem}{0}
\setcounter{lemma}{0}
\setcounter{remark}{0}
\setcounter{corollary}{0}
\setcounter{definition}{0}
\setcounter{example}{0}
}{}

\usepackage[nooneline]{caption}
\usepackage{amsmath}
\usepackage[T2A, T1]{fontenc}
\usepackage[utf8]{inputenc}
\usepackage[russian,english]{babel}
\usepackage{amsthm}
\usepackage{graphicx}
\usepackage{textcomp}
\usepackage{amsfonts}
\usepackage{amssymb}
\usepackage{epsfig}
\usepackage{color}
\usepackage{multirow}
\usepackage{dsfont}

\newtheorem{theorem}{Theorem}
\newtheorem{proposition}{Proposition}
\newtheorem{lemma}{Lemma}

\begin{document}

\makeatletter
 
\makeatother

\setcounter{page}{0}

\begin{center}
{\Large\bf On Jacobian group of the $\Delta$-graph}
\end{center}

\begin{center}
{\sc A. Mednykh$^1, ^2$, I.Mednykh$^1, ^2$, I. Yudin$^1, ^3$} \\
{\it $^1$Sobolev Institute of Mathematics, Novosibirsk, Russia\\
     $^2$Novosibirsk State University, Novosibirsk, Russia\\
     $^3$Gorno-Altaisk State University, Gorno-Altaisk, Russia\\
}
e-mail: {\tt smedn@mail.ru, ilyamednykh@mail.ru, yudin.vanya.99@mail.ru}
\bigskip

{\bf  Dedicated to the 75-th anniversary of our friend, colleague and teacher Vladislav~ Vasil'evich~Aseev}

\end{center}

\begin{abstract}
In the present paper we compute the Jacobian group of $\Delta$-graph $\Delta(n; k, l, m).$ The notion of $\Delta$-graph continues the list of families of $I$-,\, $Y$-\, and $H$-graphs well-known in the graph theory. In particular, graph $\Delta(n; 1, 1, 1)$ is isomorphic to discrete torus $C_3\times C_n.$ It this case, the structure of the Jacobian group will be find explicitly.
\smallskip


\end{abstract}

\section*{Introduction}
The notion of the Jacobian group of a graph, which is also known as the Picard group, the critical group,  the dollar or sandpile group, was independently introduced by many authors (\cite{CoriRoss}, \cite{BakerNorine}, \cite{Biggs}, \cite{BachHarpNagnib}).
This notion arises as a discrete version of the Jacobian variety from the classical theory of Riemann surfaces. It also admits a natural interpretation in various areas of physics, coding theory, and financial mathematics.

The Jacobian group is an important algebraic invariant of a finite graph. It can be defined as an   Abelian group generated by flows satisfying  the first and the second Kirchhoff lows. Equivalently, the Jacobian group can be defined as the torsion subgroup of cokernel for the Laplacian matrix of a graph. In the present paper, we will follow the latter definition. The  order of the  Jacobian group coincides with the number of spanning trees of the graph, This number  is known for some simple graphs such as the wheel, fan, prism, ladder, and M\"obius ladder \cite{BoePro}, grids \cite{NP04}, lattices \cite{SW00}, Sierpinski gaskets (\cite{AD11}, \cite{CCY07}), $3$-prism and $3$-anti-prism \cite{SWZ16}  and so on. At the same time, the structure of the Jacobian is known only in particular cases \cite{CoriRoss}, \cite{Biggs}, \cite{Lor}, \cite{YaoChinPen}, \cite{ChenHou}, \cite{MedZind}, \cite{MedMed16} and \cite{MedMed20}.

We mention that the structure of Jacobian group as well as the number of spanning trees for circulant graphs is expressed in terms of the Chebyshev polynomials. See papers  \cite{ZhangYongGol}, \cite{ZhangYongGolin}, \cite{XiebinLinZhang}, \cite{MedMed17}  and \cite{MedMed18}. More generally, this result also holds for arbitrary cyclic coverings of a graph \cite{KwonMedMed20}. In particular, (see \cite{KMM} and \cite{Ilya}) this is true for the generalized Petersen graph $GP(n, k)$ and $I$-graph $I(n,j,k).$ These two graphs are expansions of the graph consisting of a single edge. In the same time, there is a wide family of  graphs that are expansions of the tree consisting of more than one edge.  Biggs \cite{BigIYH} calls these graphs the $I$-graph, $Y$-graph, and $H$-graph, because of the method by which the graphs were created. See also \cite{HorBou} for further generalisation.

In the present paper, we investigate the structure of Jacobian for a $\Delta$-graph that is a natural counterpart of the above mentioned families of graphs. They are expansions of the triangle graph. The precise definition is as follows.

A $\Delta$-graph $\Delta(n; k,l,m)$ has $3n$ vertices $v_{x,y}$ where $x=1,2,3$ and $y$ in the integers modulo $n.$ We organize the set of vertices into three groups of $n$ vertices by first subscript $x.$ The graph induced by the vertices $v_{x,y}$ for a fixed $x,$ have edges which join each $v_{x,y}$ to $v_{x,y+j(x)},$ where $j(1)=k,\,j(2)=l$ and $j(3)=m.$ The subscript addition is performed modulo $n.$ On the other hand, the graph induced by the vertices $v_{x,y}$ for a fixed $y$ is a triangle graph $C_{3}.$ It has a $\Delta$ shape. This is a reason for a name of the graph family. This definition produces a connected graph if one has $\textrm{gcd}(k, l, m, n)=1,$ where $\textrm{gcd}(k, l, m, n)$ is the greatest common divisor of $k, l, m$ and $n.$ In what follows, we deals with connected graphs only. 

$\Delta$-graphs were investigated by the third named author in his diploma work~\cite{Yudin22}.
He found closed formulae for the number of spanning trees and the number of spanning rooted forests in such graphs.  

In this paper, we produce an approach for counting Jacobian group of the $\Delta$-graphs $  \Delta(n; k,l,m).$ The main result here is Theorem~\ref{DgraphJac} which states that group $Jac(\Delta(n;k,l,m))$ is isomorphic to the torsion subgroup  of the cokernel of a   $2\times2$ matrix whose entires are given circulant matrices.   

We note that $\Delta(n; 1,1,1)$ is isomorphic to the Cartesian product of two cyclic graphs $C_3\times C_n,$ that is can be considered as a discrete torus. It this case, the structure of the Jacobian group will be found explicitly (Theorem~\ref{jacobian}).

\section{Basic definitions and preliminary facts}

We consider a connected finite graph $G$ possibly with multiple edges but without loops. Let $V(G)$ and $E(G)$ denote the vertex set and the edge set of $G.$ For two given vertices $u, v\in V(G),$ we denote by $a_{uv}$ the number of edges between $u$ and $v.$ The matrix $A=A(G)=\{a_{uv}\}_{u, v\in V(G)}$ is called \textit{the adjacency matrix} of the graph $G.$ The degree $d(v)$ of a vertex $v \in V(G)$ is defined by $d(v)=\sum_{u\in V(G)}a_{uv}.$ Consider $D=D(G)$ diagonal matrix indexed by the elements of $V(G)$ with $d_{vv}=d(v).$ The matrix $L=L(G)=D(G)-A(G)$ is called \textit{the Laplacian matrix}, or simply \textit{Laplacian} of the graph $G.$

We refer to paper \cite{Lor} for the following helpful relation between the Laplacian matrix and the Jacobian of a graph $G.$ We consider the Laplacian $L(G)$ as a matrix of linear operator of the lattices $L(G):\mathbb{Z}^{|V|}\to\mathbb{Z}^{|V|},$ where $|V|=|V(G)|$ is the number of vertices in $G.$ The cokernel $\textrm{coker}\,(L(G))=\mathbb{Z}^{|V|}/\textrm{im}\,(L(G))$ is an Abelian group. The Smith normal form of this group is 
$$\textrm{coker}\,(L(G))\cong\mathbb{Z}_{d_{1}}\oplus\mathbb{Z}_{d_{2}}\oplus\cdots\oplus\mathbb{Z}_{d_{|V|}}.$$
It satisfies the conditions $d_i\big{|}d_{i+1},\,(1\le i\le|V|-1).$ In the case of the connected graph, the groups ${\mathbb Z}_{d_{1}},{\mathbb Z}_{d_{2}},\ldots,{\mathbb Z}_{d_{|V|-1}}$ are finite, and $\mathbb{Z}_{d_{|V|}}=\mathbb{Z}.$ We define \textit{Jacobian group} $Jac(G)$ to be the torsion subgroup of $\textrm{coker}\,(L(G)).$ In other words,
$$Jac(G)\cong\mathbb{Z}_{d_{1}}\oplus\mathbb{Z}_{d_{2}}\oplus\cdots\oplus\mathbb{Z}_{d_{|V|-1}}.$$

Let $M$ be an integer $n\times n$ matrix, then we can interpret $M$ as $\mathbb{Z}$-linear operator from $\mathbb{Z}^n$ to $\mathbb{Z}^n.$ In this interpretation $M$ has a kernel $\textrm{ker}M,$ an image $\textrm{im}\,M,$ and a cokernel $\textrm{coker} M = \mathbb{Z}^n/\textrm{im} M.$ We emphasize that $\textrm{coker}\,M$ of the matrix $M$ is completely determined by its Smith normal form. In particular, if matrices $M$ and $M^{\prime}$  are elementary equivalent then $\textrm{coker}\,M\cong\textrm{coker}\,M^{\prime}.$

In what follows, by $I_n$ we denote the identity matrix of order $n.$

We call an $n\times n$ matrix {\it circulant,} and denote it by $circ(a_0, a_1,\ldots,a_{n-1})$ if it is of the form
$$circ(a_0, a_1,\ldots, a_{n-1})=
\left(\begin{array}{ccccc}
a_0 & a_1 & a_2 & \ldots & a_{n-1} \\
a_{n-1} & a_0 & a_1 & \ldots & a_{n-2} \\
  & \vdots &   & \ddots & \vdots \\
a_1 & a_2 & a_3 & \ldots & a_0\\
\end{array}\right).$$

Recall \cite{PJDav} that the eigenvalues of matrix $C=circ(a_0,a_1,\ldots,a_{n-1})$ are given by the following simple formulas $\lambda_j=p(\varepsilon^j_n),\,j=0,1,\ldots,n-1$ where $p(x)=a_0+a_1 x+\ldots+a_{n-1}x^{n-1}$ and $\varepsilon_n$ is an order $n$ primitive root of the unity. Moreover, the circulant matrix $C=p(T),$ where $T=circ(0,1,0,\ldots,0)$ is the matrix representation of the shift operator $T:(x_0,x_1,\ldots,x_{n-2},x_{n-1})\rightarrow(x_1, x_2,\ldots,x_{n-1},x_0).$


\section{The structure of Jacobian group for the graph $\Delta(n; k,l,m)$}

We are going to investigate Jacobian group of the graph $\Delta(n; k,l,m)$ through its Laplacian matrix. Denote by $T=circ(0,1,\ldots,0)$ the $n \times n$ shift operator. Then the Laplacian $L=L(\Delta(n;  k,l,m))$ can be represented in the form
$$L=
\left(\begin{array}{ccc}
A  & -I & -I \\
-I & B  & -I  \\
 -I & -I & C \\
\end{array}\right),$$ {} 
where $A=4-T^{k}-T^{-k},\,B=4-T^{l}-T^{-l},\,C=4-T^{m}-T^{-m}$ and $I=I_n$ is the $n\times n$ identity matrix.

Recall that $Jac(\Delta(n;k,l,m))$ is given by the torsion subgroup of the $coker(L),$ where $L$ is considered as $\mathbb{Z}$-linear operator from $\mathbb{Z}^{3n}$ to itself. To find the structure of $coker(L)$ consider the following $n$-tuples of variables $x=(x_1,\ldots,x_n),\,y=(y_1,\ldots,y_n),\,z=(z_1,\ldots,z_n).$ Then, as an Abelian group,  $coker(L)$ has the following presentation
$$coker(L)=\langle x,y,z|L(x,y,z)^t=0\rangle.$$  
Hence,
$$coker(L)=\langle x,y,z|A x-y-z=0,\,-x+B y-z=0,\,-x-y+C z=0\rangle.$$
Since  $z=A x-y,$ we can eliminate the set of variables $z=(z_1,\ldots,z_n)$ from the above  presentation to get
$$coker(L)=\langle x,y|-x-A x+y+B y=0,\,-x+CA x-y-C y=0\rangle.$$

This leads to the following result.

\begin{theorem}\label{DgraphJac} The group $Jac(\Delta(n;k,l,m))$ is isomorphic to the torsion subgroup  of the $coker(M),$ where $M$ is the $2n\times2n$ block matrix given by
$$ M=\left(\begin{array}{cc}
-I-A    &  I + B \\
-I+C\,A & -I - C \\
\end{array}\right),$$
where $A=4-T^{k}-T^{-k},\,B=4-T^{l}-T^{-l}$ and $C=4-T^{m}-T^{-m}.$
\end{theorem}


\section{Explicit formulas for Jacobian of the graph $\Delta(n; 1,1,1)$}

In order to state the   results of this section we  introduce the following auxiliary function
$$\nu(n)=\begin{cases}\sqrt{7}\,U_{\frac{n}{2}-1}(\frac{5}{2}) \text{ if }n\text{ is odd,} \\ U_{\frac{n}{2}-1}(\frac{5}{2}) \text{ if }n\text{ is even.}\end{cases}$$ 
Here $U_{m-1}(x)=\frac{\sin(m\arccos(x))}{\sin{\arccos(x)}}$ is the Chebyshev polynomial of the second kind. We emphasize that $\nu(n)$ is an integer for all $n=0,\,1,\,2\,\ldots.$  

The main result of the present paper is  the following theorem.

\begin{theorem}\label{jacobian} Jacobian group $J_n=Jac(\Delta(n; 1,1,1)),\,n\ge 3$ of the graph $\Delta(n; 1,1,1)$ has the following structure
$$J_n=\mathbb{Z}_{\frac{\textrm{gcd}(n, \nu(n))}{\textrm{gcd}(n,3)}}\oplus\mathbb{Z}_{\nu(n)}^2\oplus\mathbb{Z}_{\widehat{\mu}(n)\nu(n)}\oplus\mathbb{Z}_{3\mu(n)\textrm{lcm}(n, \nu(n))},$$ where $\mu(n)$ is a $2$-periodic sequence $\{1, 7, 1, 7, \ldots\}$ and $ \widehat{\mu}(n)=\textrm{gcd}(n, 3)\mu(n).$ \end{theorem}

To find explicit formulas for Jacobian of the graph $\Delta(n; 1,1,1)$  we need to prove  a few preliminary propositions. We start with the following result which is a direct consequence of Theorem~\ref{DgraphJac}.

\begin{proposition}\label{D111Jac} The group $Jac(\Delta(n;1,1,1))$ is isomorphic to the torsion subgroup of the abelian group $coker(I+A)\oplus coker((-2I+A)(I+A)),$
where $A=4I-T-T^{-1}.$
\end{proposition}

Proof. Since  $A=B=C,$  by Theorem~\ref{DgraphJac} we have
\begin{eqnarray}
\label{cok1} coker(L)&=&\langle x,y|(-I-A)x+(I+A)y=0,\,(-I+A^2)x-(I+A)y=0\rangle \\
\nonumber &=&\langle x,y|(-I-A)x+(I+A)y=0,\,(-I+A^2)x-(I+A)x=0\rangle.\end{eqnarray}
Substituting  $u=y-x$ in (\ref{cok1}) we obtain
\begin{eqnarray}\label{cok2}
\label{cok2} coker(L)&=&\langle x,u|(I+A)u=0,\,(-2I+A)(I+A)x=0\rangle\\
\nonumber &=&coker(I+A)\oplus coker((-2I+A)(I+A)).
\end{eqnarray}

By virtue of Proposition~\ref{D111Jac}, $\textrm{coker}(L)$ splits into direct sum of two Abelian groups. We will formulate two more propositions to deal with each of them. 

\begin{proposition}\label{DIpAJac}
Let   $A=4I-T-T^{-1},$ where $I$ is $n\times n$ identity matrix and $T=circ(0,1,0,\ldots,0)$  is $n\times n$ matrix representation of the left shift operator. 
Then 
$$\textrm{coker}(I+A)=\mathbb{Z}_{\nu(n)}\oplus\mathbb{Z}_{ 3\mu(n) \nu(n) },$$ where $\mu(n)$ and $\nu(n)$ are the same as in Theorem~\ref{jacobian}.
\end{proposition}

Proof. We have $I+A=5I-T-T^{-1}=P(T),$ where $P(z)=-z^{-1}+5-z.$   Let 
$  \cal{P}=\left(\begin{array}{cc}
0   &  1\\
-1 & 5 \\
\end{array}\right) $ be the companion matrix of Laurent polynomial $P(z).$ By Lemma~1 from \cite{Ilya} cokernel of operator $I+A$ is isomorphic to corkernel of matrix ${\cal{P}}^n-I_2.$  So, $\textrm{coker}(I+A)=\mathbb{Z}_{d_1}\oplus\mathbb{Z}_{d_2/d_1},$ where $d_k $ is the greatest common divisor of $k$ by $k$ minors of the latter matrix.
By direct calculations we obtain
$${\cal{P}}^n-I_2=\left(\begin{array}{cc}
-1-U_{n-2}(5/2)   &  U_{n-1}(5/2)\\
- U_{n-1}(5/2) & -1+U_{n}(5/2) \\
\end{array}\right).$$

Since $U_{n-2}(5/2) =5U_{n-1}(5/2) -U_{n}(5/2) $ one can conclude that 
$$d_1=\textrm{gcd}(-1+U_{n}(5/2),U_{n-1}(5/2)).$$

We use the basic trigonometric identities to get $$-1+U_{n}(5/2)=2 T_{n/2+1}(5/2)U_{n/2-1}(5/2)\text{ and }U_{n-1}(5/2)=2 T_{n/2}(5/2)U_{n/2-1}(5/2).$$ 

Here,  $T_m(x)=\cos(m\arccos(x))$  is the Chebyshev polynomial of the first kind.
Then, for $n$ even we have  
\begin{eqnarray*}
\textrm{gcd}(-1+U_{n}(5/2),U_{n-1}(5/2))&=&U_{n/2-1}(5/2)\,\textrm{gcd}(2 T_{n/2+1}(5/2),2 T_{n/2}(5/2)) \\&=& \nu(n)\,\textrm{gcd}(2 T_{n/2+1}(5/2),2 T_{n/2}(5/2)).  
\end{eqnarray*}

Also, for $n$ odd we get

\begin{eqnarray*}
\textrm{gcd}(-1+U_{n}(5/2),U_{n-1}(5/2))&=&\sqrt{7}\,U_{n/2-1}(5/2)\,\textrm{gcd}(\frac{2}{\sqrt{7}} T_{n/2+1}(5/2),\frac{2}{\sqrt{7}}T_{n/2}(5/2))\\
&=&\nu(n)\textrm{gcd}(\frac{2}{\sqrt{7}} T_{n/2+1}(5/2),\frac{2}{\sqrt{7}}T_{n/2}(5/2)).
\end{eqnarray*}

Both integer sequences $2T_{k}(5/2)$ and $\frac{2}{\sqrt{7}}T_{k+1/2}(5/2))$ satisfy  the Chebyshev recursion $t_{k+1}-5 t_k +t_{k-1}=0$ with relatively prime initial data. So, we get 

$$\textrm{gcd}(2 T_{n/2+1}(5/2),2 T_{n/2}(5/2)) =1\text{ if } n \text{ is even}$$ and  
$$\textrm{gcd}(\frac{2}{\sqrt{7}} T_{n/2+1}(5/2),\frac{2}{\sqrt{7}}T_{n/2}(5/2))=1\text{ if } n \text{ is odd}.$$ 
 
 Hence, $$\textrm{gcd}(-1+U_{n}(5/2),\,U_{n-1}(5/2))=\nu(n) \text{ for all  }  n.$$
    
Also $\det({\cal{P}}^n-I_2)=-21(U_{\frac{n}{2}-1}(5/2))^2=- 3\mu(n) \nu(n)^2,$ that is $d_2=    3\mu(n) \nu(n)^2.$ So invariant factors of matrix ${\cal{P}}^n-I_2$  are $d_1=\nu(n)$ and $d_2/d_1=3\mu(n)\nu(n)$ and the result follows.
  
\begin{proposition}\label{DA2m2IApIJac}
Let  $A=4I-T-T^{-1},$ where $I$ is $n\times n$ identity matrix and $T$ is $n\times n$ matrix representation of the left shift operator $T:(x_0,\,x_1,\ldots,x_{n-2},\,x_{n-1})\to(x_1,\,x_2,\ldots,x_{n-1},\,x_{0}).$ Then 
$$\textrm{coker}(-2I+A)(I+A)=\mathbb{Z}_{\frac{\textrm{gcd}(n, \nu(n))}{\textrm{gcd}(n, 3)}}\oplus\mathbb{Z}_{\nu(n)}\oplus\mathbb{Z}_{\frac{n\widehat{\mu}(n)\nu(n)}{\textrm{gcd}(n, \nu(n))}}\oplus\mathbb{Z},$$ where $\widehat{\mu}(n)$ and $\nu(n)$ are the same as in Theorem~\ref{jacobian}.
\end{proposition}
Proof. We note that $(-2I+A)(I+A)=Q(T),$ where $Q(z)=z^{-2}-7z^{-1}+12-7z+z^2.$ Let 
$\cal{Q}=\left(\begin{array}{cccc}
0  & 1 &  0  &  0\\
0  & 0 &  1  &  0 \\
0  & 0 &  0  &  1 \\
-1 & 7 & -12 &  7 \\
\end{array}\right)$ be the companion matrix of Laurent polynomial $Q(z).$ By Lemma~1 from \cite{Ilya} cokernel of operator $Q(T)$ is isomorphic to corkernel of matrix ${\cal{B}}(n)={\cal{Q}}^n-I_4.$ Hence, $\textrm{coker}\,Q(T)=\mathbb{Z}_{d_1}\oplus\mathbb{Z}_{d_2/d_1}\oplus\mathbb{Z}_{d_3/d_2}\oplus\mathbb{Z}_{d_4/d_3},$ where  $d_k$ is the greatest common divisor of $k$ by $k$ minors of the latter matrix. First of all, we note that $d_4=\det({\cal{B}}(n))=0.$ Hence, the last term $\mathbb{Z}_{d_4/d_3}$ is equal to $\mathbb{Z}.$ So, we have to calculate the numbers $d_1,\,d_2,\,d_3$ only.

By the standard properties of companion matrices all the entries of matrix ${\cal{B}}(n)$  satisfy the following difference equation
$$x(n-2)-7x(n-1)+12x(n)-7x(n+1)+x(n+2)=0.$$ 
This equation has the   four linear independent solutions $1,\,n,\, T_n(5/2),\,U_{n-1}(5/2)$ which form a basis in the set of all solutions. We set $u=-2+2 T_n(5/2),\,v=U_{n-1}(5/2)$ and show that all the entries of ${\cal{B}}(n)$ are linear combinations of three functions $n,\,u$ and $v.$ Indeed, by direct calculations ${\cal{B}}(n)$ is given by

$$\left(\begin{array}{cccc}
\frac{n}{3} + \frac{5u}{6} - \frac{23v}{6} & -2n - \frac{11u}{6} + \frac{17v}{2} & 
2n + \frac{7u}{6} - \frac{11v}{2} & -\frac{n}{3} - \frac{u}{6} + \frac{5v}{6}\\
\frac{n}{3} + \frac{u}{6} -  \frac{5 v}{6}  &-2 n - \frac{u}{3} + 2 v  & 2 n + \frac{u}{6} - \frac{3 v}{2} &  - \frac{n}{3}+ \frac{v}{3} \\
\frac{n}{3} - \frac{v}{3} & -2 n + \frac{u}{6} +  \frac{3 v}{2} &  2 n - \frac{u}{3} - 2 v &  - \frac{n}{3} + \frac{u}{6} +  \frac{5 v}{6} \\
\frac{n}{3} - \frac{u}{6} -  \frac{5 v}{6} & -2 n +  \frac{7 u}{6} +  \frac{11 v}{2} &  2 n -  \frac{11 u}{6} - \frac{17 v}{2} & - \frac{n}{3} +  \frac{5 u}{6} +  \frac{23 v}{6} \\
\end{array}\right).$$
\bigskip

We are going to prove that $d_1=\textrm{gcd}(n,u,v)/\textrm{gcd}(n,3).$ Consider the $\mathbb{Z}$-linear span $J$ of all entries of ${\cal{B}}={\cal{B}}(n).$ It is obvious that $J$ is a main ideal in $\mathbb{Z}$ generated by $d_1.$ That is $J=d_1\mathbb{Z}.$ Consider the following three numbers $a=\frac{-u+v}{2},b=\frac{n-v}{3}$ and $c=\frac{u}{3}.$ They all belong to $J.$ Indeed $a={\cal{B}}_{1,4} + {\cal{B}}_{3, 1} - {\cal{B}}_{2, 1} - {\cal{B}}_{3, 4},\,b={\cal{B}}_{3, 1}$ and $c={\cal{B}}_{2, 1} + {\cal{B}}_{3, 4}.$ Here ${\cal{B}}_{i,j}$ is the $(i,j)$-entry of matrix ${\cal{B}}.$ Moreover, they form a basis in $J.$ To show this we note that $n=2a + 3b + 3c, u=3c$ and $v=2a + 3c.$ Also, 
$${\cal{B}}(n)=\left(\begin{array}{cccc}
-7 a + b - 8 c & 13 a - 6 b + 14 c & -7 a + 6 b - 7 c & a - b + c\\
-a + b - c & -6 b - c & a + 6 b + 2 c & -b\\
b & -a - 6 b - c & 6 b - c & a - b + 2 c\\ 
-a + b - 2 c & 7 a - 6 b + 14 c & -13 a + 6 b - 25 c & 7 a - b + 13 c
\end{array}\right).$$ So, each entry of ${\cal{B}}(n)$ is an integer linear combination of $a,b$ and $c.$ Hence $d_1=\textrm{gcd}(a,b,c)=\textrm{gcd}(\frac{-u+v}{2},\,\frac{n-v}{3},\,\frac{u}{3}).$ Analyzing divisibility of $-u+v,\,n-v,\,u$ by $2$ and $3$ we conclude that 
\begin{equation}\label{divan}
d_1=\frac{\textrm{gcd}(-u+v,\,n-v,\,u)}{\textrm{gcd}(n,3)}=\frac{\textrm{gcd}(n,u,v)}{\textrm{gcd}(n,3)}.
\end{equation}

To calculate $d_2$ we introduce matrix $\mathcal{M}=Minors({\cal{B}},2)$ consisting of all the $2$ by $2$ minors of matrix ${\cal{B}}.$ A priori, the entries of this matrix are order two polynomials in variables $n,\,u,\,v.$ This makes the calculations a bit confusing. The  following lemma helps to bypass the problem.
\begin{lemma}\label{Pythagoras} We have $u^2-21v^2=-4u.$
\end{lemma}
Proof. By definition $T_n(z)=\cos(n\arccos(z))$ and $U_{n-1}(z)=\displaystyle{\frac{\sin(n\arccos(z))}{\sin(\arccos(z))}}.$ Then, from elementary identities $$\cos(n\arccos(z))^2+\sin(n\arccos(z))^2=1\text{ and }\sin(\arccos(z))^2=1-z^2$$ we get
$$T_n(z)^2+(1-z^2)U_{n-1}(z)=1.$$ In particular case $z=5/2$ it gives  
$4T_n(5/2)^2-21U_{n-1}(5/2)^2=4.$ Note that $4T_n(5/2)^2=(u+2)^2$ and $U_{n-1}(5/2)^2=v^2.$  Hence $u^2-21v^2=-4u.$
\bigskip

Put $\displaystyle{x=-\frac{u}{18}+\frac{n\,v}{6}}$ and $\displaystyle{y=\frac{n\,v}{2}.}$ By making use of Lemma~\ref{Pythagoras} and routine calculations we obtain the following equality 
\begin{equation}\label{comb}\mathcal{M}=x\,A +y\,B,\end{equation} 
where $A$ and $B$ are given by the following integer matrices 

$$\left(\begin{array}{cccccc}
 2 + 15 n& -2 (2 + 9 n) &  2 + 3 n&  6 (1 + 3 n)& -4 - 3 n& 2 \\
 2(5 + 9 n)& -20 - 21 n&  10 + 3 n& 6 (5 + 3 n)& -20&  10 - 3 n\\
 6 (8 + 3 n)& -6 (16 + 3 n)&  48&  144&  6 (-16 + 3 n)& -6 (-8 + 3 n) \\
 2 + 3 n& -4 - 3 n&  2&  6& -4 + 3 n&  2 - 3 n\\
 10 + 3 n& -20&  10 - 3 n& -6 (-5 + 3 n)& -20 + 21 n& -2 (-5 + 9 n)\\
2& -4 + 3 n& 2 - 3 n& -6 (-1 + 3 n)& 2 (-2 + 9 n)& 2 - 15 n\\
\end{array}\right)$$
 
and

$$\left(\begin{array}{cccccc}
7 - 5 n&  2 (-4 + 3 n)&  1 - n&  -2 (-4 + 3 n)&  -1 + n&  
  0 \\
-6 (-1 + n)&  -5 + 7 n&  -1 - n&  -2 (-2 + 3 n)&  
  2&  -1 + n\\
-6 (1 + n)& 6 (3 + n)& -12& -24& -6 (-3 + n)& 
  6 (-1 + n) \\
1 - n& -1 + n& 0& 2& -1 - n& 1 + n\\
-1 - n& 
  2& -1 + n& 2 (2 + 3 n)& -5 - 7 n& 6 (1 + n)\\
0& -1 - n& 1 + n& 
  2 (4 + 3 n)& -2 (4 + 3 n)& 7 + 5 n\\
\end{array}\right) $$ respectively.
\bigskip

Now $d_2$ is the greatest common divisor of the entries of matrix $\mathcal{M}.$ Let $J_2$ be  $\mathbb{Z}$-module generated by integer linear combinations of all the entries of matrix $\mathcal{M}.$ We note that $\mathcal{M}_{1,6}=2x$ and $10 \mathcal{M}_{1, 6} + \mathcal{M}_{5, 2}= 2y.$ So, both $2x$ and $2y$ belong to $J_2.$ Also, it should be noted that $x$ and $y$ are not necessary elements of $J_2.$  
Nevertheless, if $n$ is even then   
$$-\frac{n}{2}\mathcal{M}_{5,2} + \mathcal{M}_{4,6} - \frac{7n+2}{2}\mathcal{M}_{1,6}=y.$$ That is $2x,\,y\in J_2.$ Moreover, in this case, all the entries of matrix $A$ are even integers and
\begin{equation}\label{comb}\mathcal{M}=2x\cdot\frac{1}{2}A+y\,B,\end{equation} 
The latter means that $2x$ and $y$ form a basis in $J_2$  and $d_2=\textrm{gcd}(2x,\,y).$

In a similar way, if $n$ is odd then
$$\frac{7n+11}{2}\mathcal{M}_{1,6} + \mathcal{M}_{4,6} - \frac{n+1}{2}\mathcal{M}_{5,2}=x$$
and all the entries of matrix $B$ are even integers. Now

\begin{equation}\label{comb}\mathcal{M}= x\,A +2y\cdot \frac{1}{2} B.\end{equation} 
That is $x$ and $2y$ form a basis in $J_2$ for odd $n$ and $d_2=\textrm{gcd}(x,\,2y).$

Carefully analyzing divisibility of $u$ and $n v$ by $2$ and $3$ and taking into account that  $x=-\frac{u}{18}+\frac{n\,v}{6}$ and $y=\frac{n\,v}{2}$ we finally obtain
 
\begin{eqnarray}\label{dd2}
\nonumber d_2&=&\frac{\textrm{gcd}(u,\,n v)}{\textrm{gcd}(n,3)},\,\,\,\,\,\,\text{ if }n\,(\textrm{mod}\ 3)\neq0 \\
&\text{ and } &\\
\nonumber d_2&=&\frac{1}{2}\frac{\textrm{gcd}(u,\,n v)}{\textrm{gcd}(n,3)},\,\text{ if }n\,(\textrm{mod}\ 3)=0.
\end{eqnarray}

In order to find $d_3$ we introduce matrix $\mathcal{N}=Minors[{{\cal{B}}(n)},3]$ formed by all the $3$ by $3$ minors of ${\cal{B}}(n).$ Direct calculations lead to following identity 
$$\mathcal{N}=\frac{n(u^2-21v^2)}{12}\left(\begin{array}{cccc}
1& -1& 1& -1 \\
6& -6& 6& -6 \\
6& -6& 6& -6 \\
1& -1& 1& -1 
\end{array}\right)=\frac{n\,u}{3}\left(\begin{array}{cccc}
-1& 1& -1& 1 \\
-6& 6& -6& 6 \\
-6& 6& -6& 6 \\
-1& 1& -1& 1 
\end{array}\right).$$ 
Hence, $d_3=\frac{n\,u}{3}.$
 
To finish the proof of Proposition we have to show that 
$$d_1=\frac{\textrm{gcd}(n,\nu(n))}{\textrm{gcd}(n,3)},\,d_2/d_1=\nu(n),d_3/d_2=\frac{n\widehat{\mu}(n)\nu(n)}{\textrm{gcd}(n, \nu(n))},$$ 
where $\widehat{\mu}(n)$ and $\nu(n)$ are the same as in Theorem~\ref{jacobian}.

First of all, using the arguments similar to the proof of Proposition~\ref{DIpAJac} we have $\textrm{gcd}(u,v)=\nu(n)$ if $n\,(\textrm{mod}\ 3)\neq0$ and $\textrm{gcd}(u,v)=2\nu(n)$ if $n\,(\textrm{mod}\ 3)=0.$ Also we note that  $u=3\mu\mu(n)\nu(n)^2.$ Then $u/\textrm{gcd}(u,v)=3\mu(n)\nu(n)$ if $n\,(\textrm{mod}\ 3)\neq0$ and $u/\textrm{gcd}(u,v)=\frac{3}{2}\mu(n)\nu(n)$ if $n\,(\textrm{mod}\ 3)=0.$ 

Denote by $(n)_p$ the maximum integer $\alpha$ such that $p^\alpha$ divides $n.$ By induction, one can check that $(\nu(n))_2\geq(n)_2+1$ and $(\nu(n))_3=(n)_3.$ Moreover, $(\nu(n))_7=(n)_7$ if $n$ is even and $(\nu(n))_7\leq(n)_7$ if $n$ is odd. Hence, we have $\textrm{gcd}(3\mu(n)\nu(n),n)=\textrm{gcd}(21\nu(n),n)=\textrm{gcd}(\nu(n),n)$ if $n$ is even and $\textrm{gcd}(3\mu(n)\nu(n),n)=\textrm{gcd}(3\nu(n),n)=\textrm{gcd}(\nu(n),n)$ if $n$ is odd and $n\,(\textrm{mod}\ 3)\neq0.$ 

By property (\ref{divan}) we have
$$d_1= \frac{\textrm{gcd}(n,u,v)}{\textrm{gcd}(n,3)}=\frac{\textrm{gcd}(n,\textrm{gcd}(u,v))}{\textrm{gcd}(n,3)}=\frac{\textrm{gcd}(n,\nu(n))}{\textrm{gcd}(n,3)} $$  if  $n\,(\textrm{mod}\ 3)\neq0$  and 
 $$d_1=\frac{\textrm{gcd}(n,u,v)}{\textrm{gcd}(n,3)}=\frac{\textrm{gcd}(n,\textrm{gcd}(u,v))}{\textrm{gcd}(n,3)}=\frac{\textrm{gcd}(n,2\nu(n))}{\textrm{gcd}(n,3)}=\frac{\textrm{gcd}(n,\nu(n))}{\textrm{gcd}(n,3)} $$  if  $n\,(\textrm{mod}\ 3)=0.$  

Also, the following lemma holds.
\bigskip
 
\begin{lemma}\label{gcdlemma} We have
$$\frac{\textrm{gcd}(u,n v)}{\textrm{gcd}(n,u,v)}=\textrm{gcd}(u,v).$$
\end{lemma}
 
Proof. We note that
$$\frac{\textrm{gcd}(u,\,n v)}{\textrm{gcd}(n,u,v)} =\frac{\textrm{gcd}(u,v)\textrm{gcd}(u/\textrm{gcd}(u,v),\,n v/\textrm{gcd}(u,v))}{\textrm{gcd}(n,u,v)}=\frac{\textrm{gcd}(u,v)\textrm{gcd}(u/\textrm{gcd}(u,v),\,n)}{\textrm{gcd}(\textrm{gcd}(u,v),n)}.$$

Now, we need to proof that $\textrm{gcd}(u/\textrm{gcd}(u,v),\,n)=\textrm{gcd}(\textrm{gcd}(u,v),n).$



Let $n$ be even. Note that $\mu(n)=7.$ We have $(\textrm{gcd}(u,v))_7=(\nu(n))_7=(n)_7.$ Then in the case $n\,(\textrm{mod}\ 3)\neq0 $

$$\textrm{gcd}(u/\textrm{gcd}(u,v),\,n)=\textrm{gcd}(3\cdot7\cdot\textrm{gcd}(u,v),\,n)=\textrm{gcd}(\textrm{gcd}(u,v),\,n).$$

On the other hand, if $n\,(\textrm{mod}\ 3)=0$ we have $(\textrm{gcd}(u,v))_3=(\nu(n))_3=(n)_3$ and $(\textrm{gcd}(u,v))_2=(2\nu(n))_2\geq(n)_2+2.$ Hence

$$\textrm{gcd}(u/\textrm{gcd}(u,v),\,n)=\textrm{gcd}(\frac{3\cdot7\cdot\textrm{gcd}(u,v)}{4},\,n)=\textrm{gcd}(\textrm{gcd}(u,v),\,n).$$

Now, let $n$ be odd. Then $\mu(n)=1.$ If $n\,(\textrm{mod}\ 3)\neq0$ then

$$\textrm{gcd}(u/\textrm{gcd}(u,v),\,n)=\textrm{gcd}(3\cdot\textrm{gcd}(u,v),\,n)=\textrm{gcd}(\textrm{gcd}(u,v),\,n).$$

Otherwise, if $n\,(\textrm{mod}\ 3)=0$ then, like before $(\textrm{gcd}(u,v))_3=(\nu(n))_3=(n)_3$ and $(\textrm{gcd}(u,v))_2=(2\nu(n))_2\geq(n)_2+2.$

$$\textrm{gcd}(u/\textrm{gcd}(u,v),\,n)=\textrm{gcd}(\frac{3}{4}\cdot\textrm{gcd}(u,v),\,n)=\textrm{gcd}(\textrm{gcd}(u,v),\,n).$$
   
The lemma is proved.

\bigskip
 
By making use of Lemma~\ref{gcdlemma} and equations (\ref{divan}) and (\ref{dd2}) we obtain
 
$$d_2/d_1=\frac{\textrm{gcd}(u,n v)}{\textrm{gcd}(n,u,v)}=\textrm{gcd}(u,v)=\nu(n)$$
if $n\,(\textrm{mod}\ 3)\neq0$ and
$$d_2/d_1=\frac12\frac{\textrm{gcd}(u,n v)}{\textrm{gcd}(n,u,v)}=\frac12\textrm{gcd}(u,v)=\nu(n)$$
if $n\,(\textrm{mod}\ 3)=0.$ 

Now we are going to proof that $$d_3/d_2=\frac{n\mu(n)\nu(n)}{\textrm{gcd}(n, \nu(n))}=\frac{3n\mu\mu(n)\nu(n)}{\textrm{gcd}(n,\,\nu(n))}.$$

By Lemma~\ref{gcdlemma} we have $$\frac{\textrm{gcd}(u,n v)}{\textrm{gcd}(n,u,v)}=\textrm{gcd}(u,v).$$

Therefore
$$d_3/d_2=\frac{n\,u\,\textrm{gcd}(n,3)}{3\,\textrm{gcd}(u,\,n v)}=\frac{n\,u\,\textrm{gcd}(n,3)}{3\,\textrm{gcd}(n,u,v)\,\textrm{gcd}(u,v)}=\frac{n\,u\,\textrm{gcd}(n,3)}{3\,\textrm{gcd}(n,\textrm{gcd}(u,v))\,\textrm{gcd}(u,v)}.$$

In the case $n\,(\textrm{mod}\ 3)\neq0$ we have $\textrm{gcd}(u,v)=\nu(n).$ This leads us to the following equation. 

$$d_3/d_2=\frac{n\,3\mu(n)\,\nu^2(n)\textrm{gcd}(n,3)}{3\,\textrm{gcd}(n,\nu(n))\,\nu(n)}=
\frac{n\,\mu(n)\textrm{gcd}(n,3)\,\nu(n)}{\textrm{gcd}(n,\nu(n))}=\frac{n\,\widehat{\mu}(n)\,\nu(n)}{\textrm{gcd}(n,\nu(n))}.$$

If $n\,(\textrm{mod}\ 3)=0$ then $\textrm{gcd}(u,v)=2\nu(n).$ So, we get 
$$d_3/d_2=\frac{n\,3\mu(n)\nu^2(n)\,2\,\textrm{gcd}(n,3)}{3\,\textrm{gcd}(n,2\nu(n))\,2\,\nu(n)}=\frac{n\,\mu(n)\textrm{gcd}(n,3)\,\nu(n)}{\textrm{gcd}(n,\nu(n))}=\frac{n\,\widehat{\mu}(n)\,\nu(n)}{\textrm{gcd}(n,\nu(n))}.$$

\textit{Proof} of Theorem~\ref{jacobian}.  
Combining together Propositons~\ref{D111Jac}, \ref{DIpAJac}  and \ref{DA2m2IApIJac} we already have   

$$\textrm{coker}(L)=\mathbb{Z}_{\nu(n)}\oplus\mathbb{Z}_{\frac{3\widehat{\mu}(n) \nu(n)}{\textrm{gcd}(n,3)}}\oplus 
\mathbb{Z}_{\frac{\textrm{gcd}(n, \nu(n))}{\textrm{gcd}(n, 3)}}\oplus\mathbb{Z}_{\nu(n)}\oplus\mathbb{Z}_{\frac{n\widehat{\mu}(n)\nu(n)}{\textrm{gcd}(n, \nu(n))}}\oplus\mathbb{Z}.$$

One can easily check that the following divisibility holds 
$\frac{\textrm{gcd}(n, \nu(n))}{\textrm{gcd}(n, 3)}\,\big|\,\nu(n)\,\big|\,\frac{3\widehat{\mu}(n)\nu(n)}{\textrm{gcd}(n,3)}.$
So, we need properly arrange the direct summands $\mathbb{Z}_{\frac{3\widehat{\mu}(n) \nu(n)}{\textrm{gcd}(n,3)}}$ and $\mathbb{Z}_{\frac{n\widehat{\mu}(n)\nu(n)}{\textrm{gcd}(n, \nu(n))}}.$ Using the identity
$\mathbb{Z}_{a}\oplus\mathbb{Z}_{b}=\mathbb{Z}_{\textrm{gcd}(a,b)}\oplus\mathbb{Z}_{\textrm{lcm}(a,b)}$ we have
$$\textrm{gcd}\left(\frac{3\widehat{\mu}(n) \nu(n)}{\textrm{gcd}(n,3)},\frac{n\widehat{\mu}(n)\nu(n)}{\textrm{gcd}(n, \nu(n))}\right)=\widehat{\mu}(n)\nu(n)\textrm{gcd}\left(\frac{3}{\textrm{gcd}(n,3)},\frac{n}{\textrm{gcd}(n,\nu(n))}\right)=\widehat{\mu}(n)\nu(n).$$

The last equality is true since $\textrm{gcd}(\frac{3}{\textrm{gcd}(n,3)},n)=1.$
Also, it is possible to find $\textrm{lcm}\left(\frac{3\widehat{\mu}(n) \nu(n)}{\textrm{gcd}(n,3)},\frac{n\widehat{\mu}(n)\nu(n)}{\textrm{gcd}(n,\nu(n))}\right)$ using the relation $\textrm{gcd}(a,b)\textrm{lcm}(a,b)=a\cdot b.$ From here we get 
$$\textrm{lcm}\left(\frac{3\widehat{\mu}(n) \nu(n)}{\textrm{gcd}(n,3)},\frac{n\widehat{\mu}(n)\nu(n)}{\textrm{gcd}(n,\nu(n))}\right)=\frac{3\widehat{\mu}(n)}{\textrm{gcd}(n,3)} \textrm{lcm}(n,\nu(n))=3\mu(n)\textrm{lcm}(n,\nu(n)).$$

So, in conclusion we obtain
$$\textrm{coker}(L)= \mathbb{Z}_{\frac{\textrm{gcd}(n, \nu(n))}{\textrm{gcd}(n,3)}}\oplus\mathbb{Z}_{\nu(n)}^2\oplus\mathbb{Z}_{\widehat{\mu}(n)\nu(n)}\oplus\mathbb{Z}_{3\mu(n)\textrm{lcm}(n, \nu(n))}\oplus\mathbb{Z},$$ 
where $\mu(n)$ is a $2$-periodic sequence $\{1, 7, 1, 7, \ldots\}$ and $ \widehat{\mu}(n)=\textrm{gcd}(n, 3)\mu(n).$

Since $J_n$ is given by the torsion part of $\textrm{coker}(L),$ the statement of Theorem~\ref{jacobian} follows.
 
\section*{ACKNOWLEDGMENTS}
The authors were supported by the state contract of the Sobolev Institute of Mathematics (project FWNF 2020--0005).


\begin{resetCounters}\end{resetCounters}

\end{document}